\documentclass[11pt,reqno]{article}
\usepackage{amsmath,mathrsfs,graphicx,amssymb,amsfonts,color}

\font\bg=cmbx10 scaled\magstep1

\font\small=cmr8

\newtheorem{newlemma}{{\bf Lemma}}

\newenvironment{lema}{\begin{newlemma}{\hspace{-0.5
em}{\bf.}}}{\end{newlemma}}

\newtheorem{newteorem}{{\bf Theorem}}

\newenvironment{teorem}{\begin{newteorem}{\hspace{-0.5
em}{\bf.}}}{\end{newteorem}}

\newtheorem{newkorolari}{{\bf Corollary}}

\newenvironment{korolari}{\begin{newkorolari}{\hspace{-0.5
em}{\bf.}}}{\end{newkorolari}}

\newtheorem{newdefine}{{\bf Definition}}

\newtheorem{newquestion}{{\bf Question}}

\newtheorem{newkonjek}{{\bf Conjecture}}

\newtheorem{newexample}{{\bf Example}}

\begin{document}
\tolerance=10000
\baselineskip18truept
\newbox\thebox
\global\setbox\thebox=\vbox to 0.2truecm{\hsize 0.15truecm\noindent\hfill}
\def\boxit#1{\vbox{\hrule\hbox{\vrule\kern0pt
\vbox{\kern0pt#1\kern0pt}\kern0pt\vrule}\hrule}}
\def\qed{\lower0.1cm\hbox{\noindent \boxit{\copy\thebox}}\bigskip}
\def\ss{\smallskip}
\def\ms{\medskip}
\def\bs{\bigskip}
\def\c{\centerline}
\def\nt{\noindent}
\def\ul{\underline}
\def\ol{\overline}
\def\lc{\lceil}
\def\rc{\rceil}
\def\lf{\lfloor}
\def\rf{\rfloor}
\def\ov{\over}
\def\t{\tau}
\def\th{\theta}
\def\k{\kappa}
\def\l{\lambda}
\def\L{\Lambda}
\def\g{\gamma}
\def\d{\delta}
\def\D{\Delta}
\def\e{\epsilon}
\def\lg{\langle}
\def\rg{\rangle}
\def\p{\prime}
\def\sg{\sigma}
\def\ch{\choose}

\newcommand{\ben}{\begin{enumerate}}
\newcommand{\een}{\end{enumerate}}
\newcommand{\bit}{\begin{itemize}}
\newcommand{\eit}{\end{itemize}}
\newcommand{\bea}{\begin{eqnarray*}}
\newcommand{\eea}{\end{eqnarray*}}
\newcommand{\bear}{\begin{eqnarray}}
\newcommand{\eear}{\end{eqnarray}}

\centerline{\Large \bf  On the $n-$dominating graph of specific graphs}
\bigskip

\bs

\baselineskip12truept
\centerline{Saeid Alikhani$^{}${}\footnote{\baselineskip12truept\it\small
Corresponding author. E-mail: alikhani@yazd.ac.ir}, and Davood Fatehi }
\baselineskip20truept
\centerline{\it Department of Mathematics, Yazd University}
\vskip-8truept
\centerline{\it  89195-741, Yazd, Iran}

\vskip-0.2truecm
\nt\rule{16cm}{0.1mm}

\nt{\bg ABSTRACT}
\medskip

\baselineskip14truept

\nt{ Let $G=(V,E)$ be a graph. A set $S\subseteq V(G)$ is a  dominating set, if
every vertex in $V(G)\backslash S$ is adjacent to at least one vertex in $S$. The $k$-dominating graph of $G$, $D_k (G)$, is defined to be
the graph whose vertices correspond to the dominating sets of $G$ that have cardinality
at most $k$. Two vertices in $D_k(G)$ are adjacent if and only if the corresponding dominating
sets of $G$ differ by either adding or deleting a single vertex. In this paper we consider and study  the $n$-dominating graph of specific graphs.}

\ms

\nt{\bf Mathematics Subject Classification:} {\small 05C60, 05C69.}

\nt\rule{16cm}{0.1mm}

\baselineskip20truept

\section{Introduction}
Let $G = (V,E)$ be a simple graph. 
   For any vertex $v\in V$, the open neighborhood of $v$ is the set $N(v)=\{u\in V |uv \in E\}$
    and the closed neighborhood of $v$ is the set $N[v] = N(v)\cup {v}$. For a set $S\subseteq V$, the open neighborhood
of $S$ is $N(S)=\bigcup_{v\in S} N(v)$ and the closed neighborhood of $S$ is $N[S] = N(S)\cup S$.
 A set $S \subseteq V$ is a dominating set of $G$, if $N[S] = V$, or equivalently,
 every vertex in $V\backslash S$ is adjacent to at least one vertex in $S$.
 The {\it domination number} $\gamma(G)$ is the minimum cardinality of a dominating set in $G$.
For more study in domination theory, the reader is referred to~\cite{book}. A dominating set with cardinality $\gamma(G)$ is  called a {\it
$\gamma$-set}.

\nt  Given a graph $G$, the $k$-dominating graph of $G$, $D_k (G)$, is defined to be
the graph whose vertices correspond to the dominating sets of $G$ that have cardinality
at most $k$. Two vertices in $D_k(G)$ are adjacent if and only if the corresponding dominating
sets of $G$ differ by either adding or deleting a single vertex (\cite{Hass}).

\nt A motivation for the   study of the $k$-dominating graph  is  relationships between dominating sets. In particular, given dominating
sets $S$ and $T$, is there a sequence of dominating sets $S_0 = S_1, S_2,\ldots, S_k = T$
such that each $S_{i+1}$ is obtained from $S_i$ by deleting or adding a single vertex? (see \cite{Hass}).

\nt A reconfiguration problem asks whether (when) one feasible solution to a problem
can be transformed into another by some allowable set of moves, while maintaining
feasibility at all steps. The graph $D_k(G)$
aids in studying the reconfiguration problem for dominating sets (see \cite{Hass}).
Authors in \cite{Hass} gave conditions that ensure $D_k(G)$ is connected.
\ms
\nt As usual we denote the maximum and minimum degree of a graph $G$, by $\delta(G)$ and $\Delta(G)$, respectively.

\nt In  the next section we study the $n$-dominating graph of the complete graph $K_n$. In Section 3, we consider the $n$-dominating graphs of the path graph $P_n$ and the cycle graph $C_n$  and obtain some properties of these kind
of graphs.  Finally, we study the order of dominating graph of some graph products, in Section 4.

\section{The $n$-dominating graph of $K_n$}

\nt Let $G$ be a simple graph of order $n$. Since the order of $n$-dominating graph of $G$ is the number of all dominating sets of $G$, so counting of
 the number of dominating sets of graphs is important in the study of this subject. 
We state the following theorem:

\begin{teorem}{\rm \cite{gcom,Brouwer}}\label{br}
For every connected graph $G$, the number of dominating sets is odd.
\end{teorem}

\nt The following corollary is an immediate consequence of Theorem \ref{br}:

 \begin{korolari}
For every connected graph $G$, the order of $n$-dominating graph of $G$ is odd.
\end{korolari}

\nt Author in \cite{gcom} has given a construction showing that for each  odd number $n$ there is a connected
graph $G$ such that  the number of its dominating sets is $n$. So we have the following corollary:

\begin{korolari}
For each odd number $n$, there is a connected graph $G$, such that its $n$-dominating graph has order $n$.
\end{korolari}

\nt The following theorem is about the $n$-dominating graph of $G$.

\begin{teorem}\label{2}
\begin{enumerate}
\item[(i)]  For every graph $G$ the $n$-dominating graph  $D_n(G)$ is a bipartite graph.
\item[(ii)] For every graph $G\neq O_n$, $D_n(G)$ is not regular.
\end{enumerate}
\end{teorem}

\nt{\bf Proof}.
\begin{enumerate}
\item[(i)] Let $X$ be the set of dominating sets of $G$ with odd cardinality and $Y$ be the set of dominating sets of even cardinality. It is  clear that $X\cup Y=V(D_n(G))$ and $X\cap Y=\phi$. Suppose that $A,B\in X$,
 then $(A\backslash B)\cup (B\backslash A)$ cannot be a vertex of $D_n(G)$.  Because $|A|=|B|$ or $\big||A|-|B|\big|\geq 2$. So $AB$ is not an edge of $D_n(G)$ and with similar argument we have this for two vertices in $Y$. Therefore $D_n(G)$ is a bipartite graph with parts $X$ and $Y$.

\item[(ii)] The degree of $V(G)$ as a vertex in  $D_n(G)$  is $n$. Now consider a vertex of $D_n(G)$ which is a $\gamma$-set of $G$. It's degree is $n-\gamma<n$. So the graph $D_n(G)$ is not regular.\quad\qed
 \end{enumerate}

\begin{figure}
 \hspace{5.5cm}
 \includegraphics [width=6cm,height=6.5cm]{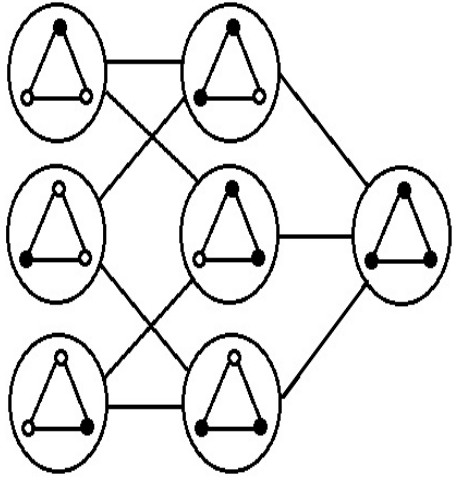}
\caption{ \label{figure2} The graph $D_{3}(K_{3})$.}
\end{figure}

\nt First we consider the $k$-dominating graph of complete graph $K_n$. It is easy to see that for every $n\in \mathbb{N}$,
$D_1(K_n)=\overline{K_n}$ and $D_2(K_n)=K_{\frac{n(n+1)}{2}}$.  Here we consider  $n-$dominating graph of $K_{n}$, $D_n(K_n)$.
The Figure \ref{figure2} shows  $D_{3}(K_{3})$.

\nt Now  we state and prove some elementary results related to
geometrical structure of $D_n(K_n)$.

 \begin{teorem}\label{pro}
 \begin{enumerate}
 \item[(i)] The number of vertices of $D_n(K_n)$ is $2^{n}-1$.
 \item[(ii)] The $n$-dominating graph $D_n(K_n)$ is a bipartite graph with parts $X$ and $Y$, with $|X|=2^{n-1}$ and $|Y|=2^{n-1}-1$.
 \item[(iii)]
 For every $n\geq 3$,  $\delta(D_{n}(K_{n}))=n-1$ and $\Delta(D_{n}(K_{n}))=n$.

 \item[(iv)]
 For every $n\geq 3$, $|E(D_n(K_n))|=n\left(2^{n-1}-1\right)$.
 \end{enumerate}
 \end{teorem}
 \nt{\bf Proof.}
 \begin{enumerate}
 \item[(i)] It is easy to see that
 the number of dominating sets of $K_n$ with cardinality $i$ is ${n \choose i}$, and since $\sum_{i=1}^n {n\choose i}=2^{n}-1$,  we have the result.

 \item[(ii)] Let $X$ be the set of dominating sets of $K_{n}$ with odd cardinality and $Y$ be the set of dominating sets of even cardinality of $K_{n}$. By theorem 2, $D_n(K_n)$ is a bipartite graph with parts $X$ and $Y$. Now we obtain the cardinality of $X$ and $Y$. Obviously $|X|=\sum_{k=1}^{\lfloor\frac{n}{2}\rfloor} {n\choose 2k-1}$ and
 $|Y|=\sum_{k=1}^{\lfloor\frac{n}{2}\rfloor} {n\choose 2k}$, therefore the result follows from well-known identities
 $\sum_{k=1}^n {n\choose k}=2^n-1$ and $\sum_{k=1}^n (-1)^k {n\choose k}=-1$.

 \item[(iii)]
For each vertex of $D_n(K_n)$, which is a dominating set of $K_n$ with cardinality one,  there exists $n-1$ dominating sets with
cardinality $2$, such that their symmetric difference is a single set, and thus both are adjacent. So, the degree of these kind of  dominating set in $D_n(K_n)$ is  $n-1$. \\
 Now suppose that $S$ is a dominating set of $K_n$ with $|S|=m$, where $2\leq m \leq n$. There are $m$ dominating sets of $K_n$ of size $m-1$ and $n-m$ dominating sets of $K_n$ of size $m+1$ which are adjacent with $S$. So, the degree of each members of $S$  in $D_n(K_n)$ is equal to $n$. Therefore $\delta(D_{n}(K_{n}))=n-1$ and $\Delta(D_{n}(K_{n}))=n$.

 \item[(iv)] By Part (i) and Hand-shaking lemma, the size of graph is $\frac{1}{2}\Big(n(n-1)+n\big(2^{n}-1-n\big)\Big)$ which is
 equal to $n\left(2^{n-1}-1\right)$.\quad\qed
\end{enumerate}

\begin{figure}[!h]
\hspace{5.5cm}
\includegraphics[width=6cm,height=6.5cm]{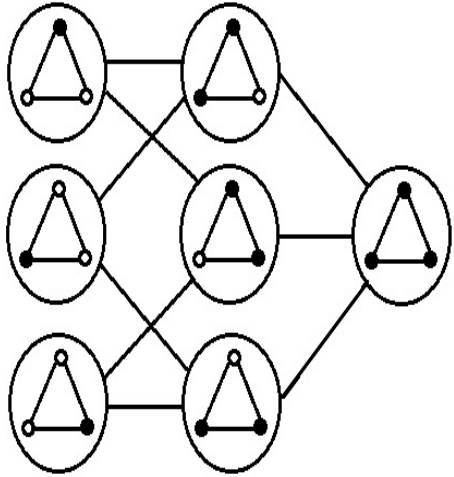}
\caption{ \label{figure1} Bipartite graph $D_{3}(K_{3})$.}
\end{figure}

\nt  We have shown the graph $D_{3}(K_{3})$ as a bipartite graph in Figure \ref{figure1}.

\nt {\bf Remark.} Since a bipartite graph with different number of vertices in its parts is not a Hamiltonian graph, so
$D_{n}(G)$ is not Hamiltonian, for any graph $G$. Also it is easy to see that for every $n\geq 3$, $D_{n}(K_{n})$ is not Eulerian and even it doesn't have an Eulerian trail.

\section{The $n$-dominating graph of path and cycle}

\nt In this section we study the $n$-dominating graph of the path and the cycle graph.

\begin{figure}[!h]
\hspace{5.5cm}
\includegraphics[width=5cm,height=.6cm]{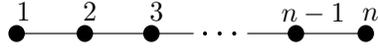}
\caption{ \label{path} The path $P_n$ with vertices labeled $\{1,2,...,n\}$.}
\end{figure}

\nt A path is a connected graph in which two vertices have degree one and the remaining vertices have degree two. Let $P_n$ be the path with
 $V (P_n) = \{1,2,...,n\}$ and $E(P_n) =\{\{1,2\},\{2,3\},...\{n-1,n\}\}$, see Figure \ref{path}.

\begin{figure}[!h]
\hspace{3.5cm}
\includegraphics[width=8.5cm,height=4.5cm]{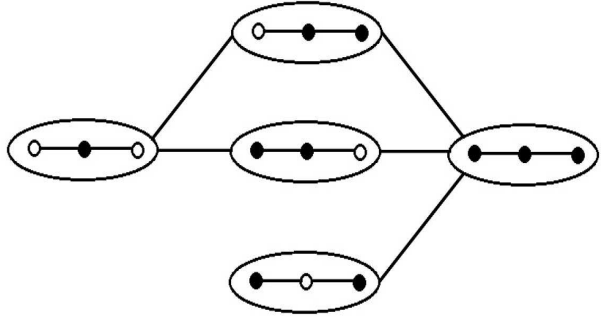}
\caption{ \label{figure3} graph $D_{3}(P_{3})$.}
\end{figure}

           \nt By the definition of  $n$-dominating graph,  $D_n(P_n)$ is a graph, which its vertices are dominating sets of the path graph $P_n$ and its two vertices are adjacent if and only if their symmetric difference is a single-member set of vertices of $P_n$. We have shown
           $D_3(P_3)$ in Figure \ref{figure3}.

\nt To study the properties of $D_n(P_n)$, we need some results about the dominating sets of $P_n$.

\begin{lema}\rm\cite{Chartrand} \label{dnp}
For every $n\in \mathbb{N}$, $\g(P_n)=\lceil\frac{n}{3}\rceil$
\end{lema}

\nt We recall that a dominating set is minimal if all its proper subsets are not dominating. Let to denote the maximum cardinality of minimal dominating sets of a graph $G$, by $\Gamma(G)$. This parameter has not been widely studied. A dominating set with cardinality $\Gamma(G)$ is  called a {\it $\Gamma$-set}. The following theorem give $\Gamma(P_n)$:

\begin{teorem}\label{G}
For every $n\in \mathbb{N}$, $\Gamma(P_n)=\lceil\frac{n}{2}\rceil$.
\end{teorem}
\nt{\bf Proof.} Consider the path $P_n$ has shown in Figure \ref{path}. It is easy to see that for odd $n$, the set $\{1,3,5,...,n\}$ is a minimal dominating set of  $P_n$ with maximum cardinality. Also for even $n$, the two sets $\{1,3,5,...,n-1\}$ and $\{2,4,6,...,n\}$ are minimal dominating sets of  $P_n$ with maximum cardinality. Obviously the cardinality of these sets  is $\lceil\frac{n}{2}\rceil$. So $\Gamma(P_n)=\lceil\frac{n}{2}\rceil$.\quad\qed

\nt By  the proof of Theorem \ref{G}, we have the following easy result:

\begin{korolari}
For  even $n\neq 4$, there are two  $\Gamma$-sets for $P_n$, and for odd $n$  there is a unique $\Gamma$-set for $P_n$.
\end{korolari}

\nt Now consider the $n$-dominating graph of $P_n$, $D_n(P_n)$. Let $\mathcal{D}(P_n,i)$ denote the family of all dominating sets of $P_n$ with cardinality $i$. Suppose that
$A,B\in \mathcal{D}(P_n,i)$ are two distinct dominating sets of $P_n$ with cardinality $i$ (or two distinct vertices of $D_n(P_n)$).  Let $\partial(A,B)$  denote the distance between  two vertices $A$ and $B$. We state and prove the following theorem:

\begin{teorem}
Suppose that
$A,B\in \mathcal{D}(P_n,i)$. We have $\partial(A,B)=2$ if and only if $\left|A\bigcap B\right|=i-1$
\end{teorem}
\nt{\bf Proof.} Let $\partial(A,B)=2$. So there exists a dominating set of $P_{n}$ with cardinality $i+1$, say $S$ such that $S\backslash A=\{x\}$ and $S\backslash B=\{y\}$. Therefore $A\bigcap B=S\backslash\{x,y\}$. Hence $$\left|A\bigcap B\right|=i+1-2=i-1.$$
Now let $\left| A\bigcap B\right|=i-1$, so there exist $x\in A$ and $y\in B$ such that $x\notin B$ and $y\notin A$. Consider $S=A\bigcup\{y\}=B\bigcup\{x\}$. The set $S$ is a dominating set of $P_n$, because $A$ and $B$ are dominating sets and hence $S$ is a vertex of $D_n(P_n)$. Note that the vertex $S$ is adjacent to $A$ and $B$, but  $A$ is not adjacent to $B$ because $\left|A\right|=\left|B\right|$. Therefore $\partial(A,B)=2$.\quad\qed

\nt Suppose that  $d(P_n,j)$, is the number of dominating sets of $P_n$ with cardinality $j$.
We need the following results:

\begin{lema}\rm\cite{Alikhani}\label{ndsp}
If $\mathfrak{p}_n^i$ is the family of dominating sets with cardinality $i$ of $P_n$, then
$$|\mathfrak{p}_n^i|=|\mathfrak{p}^{i-1}_{n-1}|+|\mathfrak{p}^{i-1}_{n-2}|+|\mathfrak{p}^{i-1}_{n-3}|.$$
\end{lema}

\begin{teorem}\rm\cite{Alikhani}\label{IJMMS}
 The following properties hold for $d(P_n,j)$:
\begin{enumerate}
\item[\textsf{(i)}]
 For every $n\in \Bbb N$, $d(P_{3n},n)=1$.

\item[\textsf{(ii)}]
For every $n\in \Bbb N$, $d(P_{3n}+2,n+1)=n+2$.

\item[\textsf{(iii)}]
For every $n\in \Bbb N$, $d(P_{3n+1},n+1)=\frac{(n+2)(n+3)}{2}-2$.

\item[\textsf{(iv)}]
For every $n\in \Bbb N$, $d(P_{3n},n+1)=\frac{n(n+1)(n+8)}{6}$.


\item[\textsf{(v)}]
If $s_n=\sum^n_{j=\lceil n/3\rceil} d(P_n,j)$, then for every $n\geq 4$, $s_n=s_{n-1}+s_{n-2}+s_{n-3}$ with initial values $s_1=1,s_2=3$ and $s_3=5$.
\item[\textsf{(vi)}]
For every $n\in \Bbb N$, $d(P_n,n-1)=n$.
\end{enumerate}
\end{teorem}

\nt We have the following properties for order of $n$-dominating graph of $P_n$ which is an immediate consequence of Theorem \ref{IJMMS}:
\begin{korolari}
For every $n\geq 4$,
$$|V(D_n(P_n)|=|V(D_{n-1}(P_{n-1})|+|V(D_{n-2}(P_{n-2})|+|V(D_{n-3}(P_{n-3})|,$$
with initial values $|V(D_1(P_1)|=1$, $|V(D_2(P_2)|=3$ and $|V(D_3(P_3)|=5$.
\end{korolari}

\nt 
 We have the following results for  the number of the dominating sets of $P_n$, i.e., $|V(D_n(P_n))|$ which are easy and straightforward to obtain:

\begin{teorem}
\begin{enumerate}
\item[(i)]
 The generating function of $|V(D_n(P_n))|$  is $T(x)=-\frac{(x+1)^2}{x^3+x^2+x-1}.$
 \item[(ii)]
 The order of $D_n(P_n)$ satisfies  $$|V(D_n(P_n))|=\frac{(-1)^n\left(A\left(\tau_1^{-1}\right)^n\tau_2\tau_3+B\left(\tau_2^{-1}\right)^n\tau_1\tau_3+
C\left(\tau_3^{-1}\right)^n\tau_1\tau_2\right)}{\tau_1\tau_2\tau_3},$$
where $\tau_1,~\tau_2$ and $\tau_3$ are the roots of $x^3+x^2+x-1=0$ and $A,~B$ and $C$ are
$\frac{\left(\tau_1-1\right)^2}{\left(\tau_3-\tau_1)(\tau_2-\tau_1\right)},~
\frac{\left(\tau_2-1\right)^2}{\left(\tau_3-\tau_2)(\tau_1-\tau_2\right)}$ and
$\frac{\left(\tau_3-1\right)^2}{\left(\tau_2-\tau_3)(\tau_1-\tau_3\right)}$, respectively.
\end{enumerate}
\end{teorem}

\nt Here we state and prove more results on the number of $\g$-sets and $(\g+1)$-sets of paths:

\begin{korolari}
The number of $\g$-sets of $P_{3k}$, $P_{3k+1}$ and $P_{3k+2}$ are $1$, $\frac{k^2+5k+2}{2}$ and $k+2$, respectively.
\end{korolari}
\nt{\bf Proof.} It follows from Lemma \ref{dnp} and Theorem \ref{IJMMS} (i),(ii) and (iii).\quad\qed

\nt We know that one or two members of the family of dominating sets with cardinality $\g+1$  are $\Gamma$-sets. We shall find the cardinality of these families. Our objective is finding any $\g$-set by adding or deleting some vertices of another $\g$-set.

\nt By Theorem \ref{IJMMS} (iv), the number of $(\g+1)$-sets of $P_{3n}$ is $d(P_{3n},n+1)=\frac{n(n+1)(n+8)}{6}$. Here we obtain the number of
$(\gamma+1)$-sets of $P_{3n+1}$ and $P_{3n+2}$:

\begin{teorem}
For every $n\in \Bbb N$,
\begin{enumerate}
\item[\textsf{(i)}]
the number of $(\g+1)$-sets of $P_{3n+1}$is $\frac{1}{24}\left(n^4+18n^3+71n^2+78n+24\right)$.
\item[\textsf{(ii)}]
the number of $(\g+1)$-sets of $P_{3n+2}$ is $\frac{1}{120}n(n+1)\left(n^3+24n^2+121n+94\right)$.
\end{enumerate}
\end{teorem}
\nt{\bf Proof.}
\begin{enumerate}
\item[\textsf{(i)}] We shall prove $d(P_{3n+2},n+2)=\frac{1}{24}\left(n^4+18n^3+71n^2+78n+24\right)$.
We use  induction on $n$. Since $d(P_5,3)=8$, the result is true for $n=1$. Now suppose that the result is true for all natural numbers less than n, and we prove it for n. By Lemma \ref{ndsp}, $$d(P_{3n+2},n+2)=d(P_{3n+1},n+1)+d(P_{3n},n+1)+d(P_{3n-1},n+1).$$
Since $d(P_{3n-1},n+1)=d(P_{3(n-1)+2},(n-1)+2)$, by Theorem \ref{IJMMS} (iii),(iv) and induction hypothesis, we have

\begin{eqnarray*}
d(P_{3n+2},n+2)&=&\frac{(n+2)(n+3)}{2}-2+\frac{n(n+1)(n+8)}{6}+\frac{1}{24}n(n^3+14n^2+23n-14)\\
&=&\frac{1}{24}\left(n^4+18n^3+71n^2+78n+24\right).
\end{eqnarray*}

\item[\textsf{(ii)}]
We use induction on $n$. Since $d(P_4,3)=4$, The result is true for $n=1$. Now suppose that the result is true for all natural numbers less than $n$, and we prove it for $n$. By Lemma \ref{ndsp}, $$d(P_{3n+1},n+2)=d(P_{3n},n+1)+d(P_{3n-1},n+1)+d(P_{3n-2},n+1).$$
Since
$d(P_{3n-1},n+1)=d(P_{3(n-1)+2},(n-1)+2)$
and
$d(P_{3n-2},n+1)=d(P_{3(n-1)+1},(n-1)+2),$
by Theorem \ref{IJMMS} (i), (iv) and induction hypothesis, we have

\begin{eqnarray*}
d(P_{3n+1},n+2)&=&\frac{n(n+1)(n+8)}{6}+\frac{1}{24}n(n^3+14n^2+23n-14)\\
&+&\frac{1}{120}(n-1)n(n^3+21n^2+76n-4)\\
&=&\frac{1}{120}n(n+1)\left(n^3+24n^2+121n+94\right).\quad\qed
\end{eqnarray*}
\end{enumerate}

\nt The following theorem is about the minimum and maximum degree of $D_n(P_n)$:

\begin{teorem}\label{delta}
For $n$-dominating graph of $P_n$ we have $\Delta(D_n(P_n))=n$ and $\delta(D_n(P_n))=n-\lceil\frac{n}{2}\rceil$.
\end{teorem}
\nt{\bf Proof.} By Theorem \ref{IJMMS}(vi), $D_n(P_n)$ has $n$ vertices of cardinality $n-1$ and these are  all  vertices adjacent to vertex $V(P_n)$, and so $\deg(V(P_n))=n$. For every $\lceil n/3\rceil<i<n$, there are at most $i$ vertices  $A\in\mathcal{D}(P_n,i-1)$ and exactly $n-i$ vertices  $B\in\mathcal{D}(P_n,i+1)$ adjacent to $C\in \mathcal{D}(P_n,i)$. For all vertices in $\mathcal{D}(P_n,\lceil n/3\rceil)$ there are $n-\g$ adjacent vertices of $\mathcal{D}(P_n,\g+1)$ and for vertices in $\mathcal{D}(P_n,\Gamma)$ there are $n-\Gamma(P_n)$ adjacent vertices of $\mathcal{D}(P_n,\Gamma+1)$. So $\delta(D_n(P_n))=n-\Gamma$ and $\Delta(D_n(P_n))=n$. We have the result by Theorem \ref{G}.\quad\qed

\begin{teorem}\label{connectedness}
The graph $D_n(P_n)$ is a connected graph.
\end{teorem}
\nt{\bf Proof.} Let $A\in \mathcal{D}(P_n,i)$ and $B\in \mathcal{D}(P_n,j)$ are two vertices of $D_n(P_n)$, where $i<j$. Step by step, we add a vertex to $A$ to reach a vertex in $\mathcal{D}(P_n,n)$. Then step by step remove a vertex to reach to a vertex in
$\mathcal{D}(P_n,j)$. Note that all the sets obtained in these steps are dominating sets of $P_n$ and so  are vertices of $D_n(P_n)$. Hence we have found a path between $A$ and $B$. Therefore $D_n(P_n)$ is a connected graph.\quad\qed

\medskip

\begin{figure}[h]\label{cycle}
\hspace{6cm}
\includegraphics[width=3.2cm, height=3.cm]{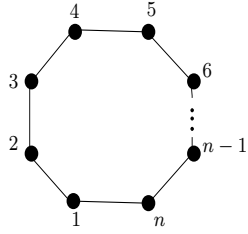}
\caption{\label{figure2'} The cycle $C_n$ with vertices labeled $\{1,2,...,n\}$.}
\end{figure}

\nt
 Here we consider cycles. Let $C_n, n\geq 3$, be the cycle with $n$ vertices $V(C_n)=\{1,2,...,n\}$ and
 $E(C_n)=\{\{1,2\},\{2,3\},\ldots,\{n-1,n\},\{n,1\}\}$, see Figure~\ref{figure2'}.

\nt Suppose that ${\cal D}(C_{n},i)$
 or simply ${\cal C}_{n}^{i}$ denote the family of dominating set of $C_{n}$ with
 cardinality $i$. The following theorem gives a recurrence relation for the number of dominating sets of cycles with cardinality $i$:

\begin{teorem}\label{theorem18}{\rm\cite{DPC}}
\begin{enumerate}
 \item[(i)] If ${\cal C}_{n}^{i}$ is the family of dominating set of $C_{n}$ with cardinality $i$, then
$|{\cal C}_{n}^{i}|=|{\cal C}_{n-1}^{i-1}|+|{\cal C}_{n-2}^{i-1}|+|{\cal C}_{n-3}^{i-1}|.$

 \item[(ii)] If $S_n=\sum_{j=\lceil\frac{n}{3}\rceil}^n |{\cal C}_n^j|$, then for every $n\geq 4$, $S_n=S_{n-1}+S_{n-2}+S_{n-3}$
with initial values $S_1=1, S_2=3$ and $S_3=7$.
\end{enumerate}
\end{teorem}

\begin{figure}[!h]\label{ddcycle}
 \hspace{5.5cm}
 \includegraphics [width=8cm,height=7.5cm]{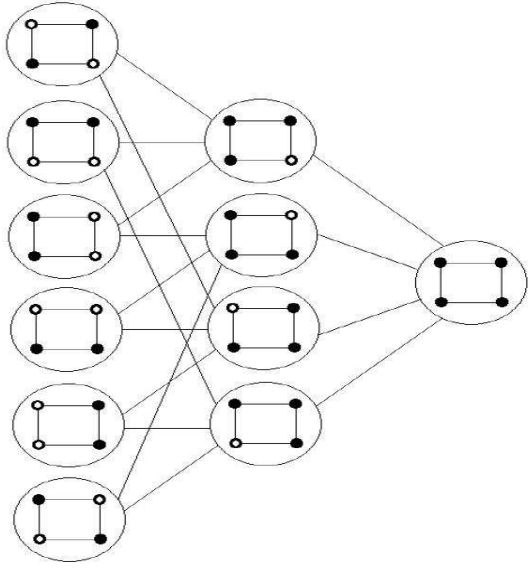}
\caption{\label{ddcycle} graph $D_{4}(C_{4})$.}
\end{figure}

\nt We have the following properties for order of $n$-dominating graph of $C_n$ which is an immediate consequence of Theorem \ref{theorem18} (see Figure \ref{ddcycle}):
\begin{korolari}
For every $n\geq 4$,
$$|V(D_n(C_n)|=|V(D_{n-1}(C_{n-1})|+|V(D_{n-2}(C_{n-2})|+|V(D_{n-3}(C_{n-3})|,$$
with initial values $|V(D_1(C_1)|=1$, $|V(D_2(C_2)|=3$ and $|V(D_3(C_3)|=5$.
\end{korolari}

\begin{teorem}
For every $n\in\mathbb N$, $\Gamma(C_n)=\lfloor\frac{n}{2}\rfloor$.
\end{teorem}
\nt{\bf Proof.} Consider the cycle $C_n$ has shown in Figure \ref{figure2'}. It is easy to see that for even $n$, the two sets $\{1,3,5,...,n-1\}$ and $\{2,4,6,...,n\}$ are minimal dominating sets of $C_n$ with maximum cardinality. For odd $n$ all the sets of the form $\{i\oplus 2k: \ 1\leq i \leq n, 1\leq k\leq \lfloor\frac{n}{2}\rfloor\}$  are the $\Gamma$-sets of $C_n$, where $\oplus$ is sum modulo $n$.  Obviously the cardinality of these sets is $\lfloor\frac{n}{2}\rfloor$. So $\Gamma(C_n)=\lfloor\frac{n}{2}\rfloor$.\quad\qed

\nt Similar to  $n$-dominating of $P_n$, we have the following results:

\begin{teorem}
\begin{enumerate}
\item[(i)]
For  even $n>4$, there are two  $\Gamma$-sets and for odd $n$ there are $n$ $\Gamma$-sets for $C_n$.
\item[(ii)]
$\Delta(D_n(C_n))=n$ and $\delta(D_n(C_n))=n-\Gamma$.
\item[(iii)]
The graph $D_n(C_n)$ is a connected graph.
\end{enumerate}
\end{teorem}

\nt The following result is easy to obtain:

\begin{teorem}
\begin{enumerate}
\item[(i)]
 The generating function of $|V(D_n(C_n))|$  is $T(x)=-\frac{3x^2+2x+1}{x^3+x^2+x-1}.$
 \item[(ii)]
 The order of $D_n(C_n)$ satisfies  $$|V(D_n(C_n))|=\frac{(-1)^n\left(A\left(\tau_1^{-1}\right)^n\tau_2\tau_3+B\left(\tau_2^{-1}\right)^n\tau_1\tau_3+
C\left(\tau_3^{-1}\right)^n\tau_1\tau_2\right)}{\tau_1\tau_2\tau_3},$$
where $\tau_1,~\tau_2$ and $\tau_3$ are the roots of $x^3+x^2+x-1=0$ and $A,~B$ and $C$ are
$\frac{3\tau_1^2-2\tau_1+1}{\left(\tau_3-\tau_1)(\tau_2-\tau_1\right)},~
\frac{3\tau_2^2-2\tau_2+1}{\left(\tau_3-\tau_2)(\tau_1-\tau_2\right)}$ and
$\frac{3\tau_3^2-2\tau_3+1}{\left(\tau_3-\tau_2)(\tau_3-\tau_1\right)}$, respectively.
\end{enumerate}
\end{teorem}

\section{The order of $n$-dominating graph of  graph products}

\nt We would like to end this paper by stating  some results on the order of $n$-dominating graph of some operations of two graphs. The join
$G = G+H$ of two graph $G$ and $H$ with disjoint vertex sets $V_1$ and $V_2$ and
edge sets $E_1$ and $E_2$ is the graph union $G\cup H$ together with all the edges joining $V_1$ and
$V_2$. The corona of two graphs $G$ and $H$ is the graph $G\circ H$ formed from one copy of $G$ and $|V(G)|$  copies of $H$,
where the $i$-th vertex of $G$ is adjacent to every vertex in the $i$-th copy of $H$.

\nt The following theorem gives the order of number of $D_n(G+H)$ and $D_n(G\circ H)$. This theorem can be found in domination polynomial format  in \cite{Oper}:

\begin{teorem}
Let $G$ and $H$ be graphs of order $p$ and $q$, respectively, and $n=p+q$. Then
\begin{enumerate}
\item[(i)] $|V(D_n(G+H))|=(2^p-1)(2^q-1)+|V(D_p(G))|+|V(D_q(H))|$.
\item[(ii)] $|V(D_n(G\circ H))|=(2^q +|V(D_q(H))|)^p$.

\end{enumerate}
\end{teorem}
\nt{\bf Proof.}
\begin{enumerate}
\item[(i)]
Let $i$ be a natural number $1\leq ≤ i \leq  p+q$. We shall count the number of dominating sets of $G+H$.
 If $i_1$ and $i_2$ are two natural numbers such that $i_1 + i_2 = i$,
then clearly, for every $D_1 \subseteq V (G)$ and $D_2 \subseteq V (H)$, such that $|D_j| = i_j$,
$j = 1, 2$, $D_1\cup D_2$ is a dominating set of $G+H$. Also, if $D$ is a dominating set of $G$ of size $i$, then $D$ is a dominating set for
 $G+H$ of size $i$. The same is true for every dominating $D$ of size $i$. Thus we have the result.

 \item[(ii)]
 In the corona of two graphs $G$ and $H$, every
vertex $u\in V(G)$ is adjacent to all vertices of the corresponding
copy of $H$. So, we can delete all edges in $E$ in the
corona. Therefore, the arising graph is the disjoint union of
$p$ copies of the corona $K_1\circ H$.  Therefore, the number of dominating sets of $G\circ H$ is $|V(D_{q+1}(K_1\circ H)|^p$. It remains
to obtain the number of dominating sets of $K_1\circ H$.
We have two cases for a dominating set $S$ of every graph of  form $K_1\circ H$.

\nt Case 1. $S$ includes the vertex originally in $K_1$ and an
arbitrary subset of the $q$ vertices from the copy of $H$. In this case we have $2^q$ dominating sets.

\nt Case 2. $S$ does not include the vertex in $K_1$ and it is exactly a dominating set
of $H$. In this case there are $|V(D_q(H))|$ dominating sets.

\nt By addition principle we have $|V(D_{q+1}(K_1\circ H)|=2^m+ |V(D_q(H))|$.\quad\qed
\end{enumerate}

\nt
Given any two graphs $G$ and $H$ we define the Cartesian product, denoted $G\Box H$,
to be the graph with vertex set $V(G)\times V (H)$ and edges between two vertices $(u_1,v_1)$
and $(u_2,v_2)$ if and only if either $u_1 = u_2$ and $v_1v_2 \in E(H)$ or $u_1u_2 \in E(G)$ and
$v_1 = v_2$. We think that computing of the order of Cartesian product of two graphs $G$ and $H$, i.e. the number of
dominating sets of $G\Box H$ is not easy. Using Theorem 4.2. of \cite{prod} we have the following recurrence relation for
 the order of Ladder graph $L_n$ which is $P_n\Box K_2$.

 \begin{teorem}{\rm\cite{prod}}
 Let $G$ be a graph of order $n$, then
 \begin{eqnarray*}
 |V(D_{2n}(L_n))|&=&3|V(D_{2n-2}(L_{n-1}))|+2|V(D_{2n-4}(L_{n-2}))|+2|V(D_{2n-6}(L_{n-3}))|\\
&-&|V(D_{2n-8}(L_{n-4}))|-|V(D_{2n-10}(L_{n-5}))|.
\end{eqnarray*}
 \end{teorem}


\nt In this section we studied the order of  $|V|$-dominating graph of some graph operations. Study the geometrical properties and parameters
of these kind of graphs  remain as open problems which we believe deserve attention.

\end{document}